\documentclass[a4paper,11pt]{article}

\usepackage{MyLaTeX,MyArticle}
\input xy
\xyoption{all}

\title{Algebraic Modules and the Auslander--Reiten Quiver}
\author{David A.~Craven, University of Oxford}
\date{December 2007}

\newcommand{\cx}{\operatorname{cx}}

\begin{document}
\maketitle

\noindent Recall that an algebraic module is a $KG$-module that satisfies a polynomial with integer coefficients, with addition and multiplication given by direct sum and tensor product. In this article we prove that non-periodic algebraic modules are very rare, and that if the complexity of an algebraic module is at least 3, then it is the only algebraic module on its component of the (stable) Auslander--Reiten quiver. We include a strong conjecture on the relationship between periodicity and algebraicity.

\section{Introduction}

Trying to decompose the tensor product of two (even simple) modules is, in general, a hopeless proposition. In some cases it might be possible to have some control over which summands can appear; following Alperin in \cite{alperin1976b}, we define a module to be \emph{algebraic} if it satisfies a polynomial with integer coefficients, where addition and multiplication are given by the direct sum and the tensor product. It is clear that a module $M$ is algebraic if and only if there are only finitely many isomorphism types of indecomposable summand in the collection of modules $M^{\otimes n}$ for all $n\geq 0$. Examples include all projective modules, more generally all trivial source modules, and all simple modules for $p$-soluble groups \cite{feit1980}.

In this article we will produce results on how the concept of algebraic modules can be related to that of the Heller operator, and how some strong results can be achieved concerning algebraic modules on the Auslander--Reiten quiver.

\begin{thma}\label{comp3} Let $\Gamma$ be a component of the stable Auslander--Reiten quiver $\Gamma_s(KG)$. Suppose that the complexity of modules on $\Gamma$ is at least 3. (In this case, $\Gamma$ has tree class $A_\infty$.) Then $\Gamma$ contains at most one algebraic module, and such a module lies on the end of $\Gamma$.
\end{thma}

\begin{thma}\label{algheller} Let $M$ be an algebraic module for a finite group $G$.
\begin{enumerate}
\item If $M$ is periodic then $\Omega^i(M)$ is algebraic for any $i\in\Z$.
\item If $M$ is non-periodic, then at most one of the modules $\Omega^i(M)$ is algebraic, and if $M$ is self-dual and one of the $\Omega^i(M)$ is algebraic, then it is $M$ that is algebraic.
\end{enumerate}
Furthermore, all possibilities allowed by this theorem do occur.
\end{thma}

This theorem broadly says that non-periodic algebraic modules are `rare', and this can be used to prove the following fairly strict condition on the number of algebraic modules in a non-projective component of the stable Auslander--Reiten quiver, which we denote by $\Gamma_s(KG)$. The next theorem appears technical, and we will single out two special cases as corollaries to the theorem.

\begin{thma}\label{vertexcpcp} Let $K$ be a field of characteristic $p$, and let $G$ be a finite group. Let $M$ be a non-periodic, indecomposable, algebraic module belonging to a wild block, and suppose that there is a subgroup $Q$, not containing a vertex of $M$, such that $M\res Q$ is non-periodic. Then no other module on the same component of $\Gamma_s(KG)$ as $M$ is algebraic.
\end{thma}

In this case, $M$ lies on a component of $\Gamma_s(KG)$ of type $A_\infty$, and we also show that $M$ lies at the end of its component.

The first corollary is Theorem \ref{comp3} itself, and the second is the following.

\begin{cora}\label{p'dim} Let $K$ be a field of characteristic $p$ and let $G$ be a finite group whose Sylow $p$-subgroups are neither tame nor isomorphic to $Q=C_p\times C_p$. Let $M$ be a non-periodic indecomposable $KG$-module whose dimension is prime to $p$. If $M$ is algebraic, then no other module on the component of $\Gamma_s(KG)$ is algebraic.
\end{cora}

In the case of the group $C_p\times C_p$, little is known. However, conjecturally there is a strong link between periodicity and whether a module is algebraic.

\begin{conja}\label{cpcpper} Let $K$ be a finite field of characteristic $p$, and let $G$ be the group $C_p\times C_p$. Let $M$ be an absolutely indecomposable module such that $p\mid \dim M$. Then $M$ is algebraic if and only if it is periodic.
\end{conja}

The reason behind the presence of a finite field is that it does not appear clear if it is merely the dimensions of indecomposable summands of powers of the module $M$ that are bounded, rather than their coming from a finite list. In the case where the field is finite, clearly both concepts coincide. We will provide our evidence for this conjecture in the final section.

The structure of this article is simple: in the following section the preliminary results needed on algebraic modules are collated. In the short succeeding section, Theorem \ref{algheller} is proved, and in the section after that we prove Theorem \ref{vertexcpcp}. The final section contains the aforementioned evidence behind Conjecture \ref{cpcpper}.

\section{Preliminaries}

In this section we will describe the preliminary results on algebraic modules, together with a result on tensor products. We start with algebraic modules, and the following lemma is easy.

\begin{lem}[{{\cite[Section II.5]{feit}}}] Let $M=M_1\oplus M_2$ be a $KG$-module, and suppose that $H_1\leq G\leq H_2$.
\begin{enumerate}
\item $M$ is algebraic if and only if $M_1$ and $M_2$ are algebraic.
\item The module $M_1\otimes M_2$ is algebraic.
\item The modules $M_1\res {H_1}$ and $M_1\ind {H_2}$ are algebraic.
\end{enumerate}
\end{lem}

An easy corollary of this lemma is that an indecomposable module is algebraic if and only if its source is.

We also need the fact that a module is algebraic if and only if it is algebraic in the stable module category.

\begin{prop}\label{quotbyprojs} Let $\ms I$ be an ideal of algebraic modules in the Green ring $a(KG)$, and let $M$ be a $KG$-module. Then $M$ is algebraic in $a(KG)$ if and only if $M+\ms I$ is algebraic in $a(KG)/\ms I$. In particular, if $\ms P$ denotes the ideal consisting of all projective modules, then a $KG$-module $M$ is algebraic if and only if $M + \mathscr P$ is algebraic.
\end{prop}
\begin{pf} Suppose that $M$ is algebraic. Then $M$ satisfies some polynomial in the Green ring, and therefore its coset in any quotient satisfies this polynomial as well. Conversely, suppose that $M+\ms I$ satisfies some polynomial in the quotient $a(KG)/\ms I$. Thus
\[ \sum \alpha_i (M+\ms I)^i=\ms I.\]
This implies that, since $(M+\ms I)^i=M^{\otimes i}+\ms I$, then
\[ \sum \alpha_i M^{\otimes i} \in \ms I,\]
which consists of algebraic modules. Hence there is some polynomial involving only $M$ witnessing the algebraicity of $M$.
\end{pf}

In fact, one can extend the ideal $\ms P$ to one containing not only the projective modules but all modules of cyclic vertex.

Since we are relating tensor products and the Heller operator, we need the next well-known lemma.

\begin{lem} Let $M$ and $N$ be modules. Then
\[ \Omega(M\otimes N)=\Omega^0(\Omega(M)\otimes N).\]
\end{lem}

We also need two results regarding summands of tensor powers, due to Benson--Carlson and Auslander--Carlson, which are necessary for the proof of Theorem \ref{algheller}. We amalgamate them into a single theorem.

\begin{thm}
Let $G$ be a finite group and let $M$ and $N$ be absolutely indecomposable $KG$-modules, where $K$ is a field of characteristic $p$.
\begin{enumerate}
\item \textbf{(\cite[Theorem 2.1]{bensoncarlson1986})} $K|M\otimes N$ if and only if $p\nmid \dim M$ and $M\cong N^*$, in which case $K\oplus K$ is not a summand of $M\otimes N$. If $p\mid \dim M$, then every summand of $M\otimes N$ has dimension a multiple of $p$. 
\item \textbf{(\cite[Proposition 4.9]{auslandercarlson1986})} If $\dim M$ is a multiple of $p$, then $M\oplus M$ is a direct summand of $M\otimes M^*\otimes M$.
\end{enumerate} 
Therefore for all $KG$-modules $M$, we have $M|M\otimes M^*\otimes M$.
\end{thm}

\section{Algebraicity and Periodicity}

In this section we will relate the Heller operator and algebraic modules. All modules are algebraic if $G$ has cyclic Sylow $p$-subgroups. If $G$ does not, then there are infinitely many non-algebraic $KG$-modules.

\begin{prop}\label{omegaknonalg} Let $G$ be a finite group of $p$-rank at least 2, and let $K$ be a field of characteristic $p$. Then, for all $i\neq 0$, the module $\Omega^i(K)$ is not algebraic.
\end{prop}
\begin{pf} If $G$ has $p$-rank 2 or more, then the trivial module, $K$, is non-periodic. Notice that, modulo projective modules,
\[ \(\Omega^i(K)\)^{\otimes n}=\Omega^{ni}(K),\]
and so $\Omega^{ni}(K)$ appears as a summand of the $n$th tensor power of $\Omega^i(K)$ for all $n\geq 1$, an infinite collection of summands since $K$ is not periodic.\end{pf}

If $G$ is not of $p$-rank 2 and does not have cyclic Sylow $p$-subgroups, then $p=2$ and the Sylow $2$-subgroups of $G$ are generalized quaternion. In this case, by the Brauer--Suzuki theorem, $G$ possesses a normal subgroup $\mathrm{Z}^*(G)$ such that $G/\mathrm{Z}^*(G)$ has dihedral Sylow $2$-subgroups, and so there are non-algebraic modules for this quotient. Alternatively, a generalized quaternion $2$-group possesses a $V_4$ quotient, and so there are non-algebraic modules for generalized quaternion $2$-groups, whence any indecomposable module for $G$ with one of those modules as a source would be non-algebraic.

Now suppose that a $KG$-module $M$ is periodic; in the next proposition, we use the obvious fact that a module $M$ is algebraic if and only if $M^{\otimes i}$ is algebraic for some $i\geq 1$.

\begin{prop} Let $M$ be an algebraic periodic module. Then $\Omega^i(M)$ is algebraic for all $i$.\end{prop}
\begin{pf} Suppose that $\Omega^n(M)=M$. We know that
\[ \Omega(M\otimes N)=\Omega^0(\Omega(M)\otimes N)=\Omega^0(M\otimes \Omega(N)).\]
Hence, $\Omega^0(\Omega^i(M)^{\otimes n})=\Omega^{ni}(M^{\otimes n})=\Omega^0(M^{\otimes n})$, and since $M^{\otimes n}$ is algebraic (as $M$ is), the module $\Omega^i(M)$ is algebraic for all $i$ (as $\Omega(M)^{\otimes n}$ is).
\end{pf}

Both possibilities allowed---that the $\Omega$-translates of $M$ are either all algebraic modules or all non-algebraic modules---occur in the module category of the quaternion group. Firstly, the trivial module is an algebraic periodic module, and secondly, since the group $V_4$ has $2$-rank 2, the non-trivial Heller translates of the trivial module for that group are non-algebraic by Proposition \ref{omegaknonalg}, and so those modules, viewed as modules for the quaternion group, are also non-algebraic. It should be mentioned that no examples of non-algebraic periodic modules are known if the characteristic of the field is odd.

Now we consider non-periodic modules. Since a module $M$ is non-periodic if and only if $M\otimes M^*$ is, we firstly consider self-dual non-periodic modules, then apply this to the general case.

\begin{prop} Let $M$ be a self-dual non-periodic module. If $i\neq 0$ then $\Omega^i(M)$ is not algebraic.\end{prop}
\begin{pf} Consider the module
\[ \Omega^0(\Omega^i(M)\otimes \Omega^i(M)\otimes \Omega^i(M))=\Omega^{3i}(M^{\otimes 3});\]
as $M$ is a summand of $M^{\otimes 3}$, we see that $\Omega^{3i}(M)$ is a summand of $\Omega^i(M^{\otimes 3})$. We can clearly iterate this procedure to prove that infinitely many different $\Omega$-translates of $M$ lie in tensor powers of $\Omega^i(M)$ (and these all contain different indecomposable summands as $M$ is non-periodic) proving that $\Omega^i(M)$ is non-algebraic, as required.
\end{pf}

\begin{cor} Let $M$ be a non-periodic algebraic module. Then no module $\Omega^i(M)$ for $i\neq 0$ is algebraic.\end{cor}
\begin{pf} Suppose that both $M$ and $\Omega^i(M)$ are algebraic. Then so is $M^*$, and therefore so is 
\[ \Omega^0(M^*\otimes\Omega^i(M))=\Omega^i(M\otimes M^*).\]
Since $M\otimes M^*$ is self-dual, this module cannot be algebraic, a contradiction.
\end{pf}

Hence for non-periodic modules $M$, either none of the modules $\Omega^i(M)$ is algebraic, or exactly one module is, and in the latter case, if one of the modules is self-dual then this is the algebraic module. In the case of the dihedral $2$-groups, there are non-periodic modules $M$ such that no $\Omega^i(M)$ are algebraic, and there are self-dual, non-periodic algebraic modules. 
This completes the proof of Theorem \ref{algheller}. This theorem has the following corollary, which is useful when computing examples.

\begin{cor}\label{easynonalg} Let $M$ be a non-periodic indecomposable module, and suppose that there is some $n\geq 2$ such that $\Omega^i(M)$ or $\Omega^i(M^*)$ is a summand of $M^{\otimes n}$ for some $i\neq 0$. Then the module $\Omega^i(M)$ is non-algebraic for all $i\in \Z$.
\end{cor}
\begin{pf}Suppose that $\Omega^i(M)$ is a summand of $M^{\otimes n}$, for some $n\geq 2$ and $i\neq 0$. Then, for each $j\in\Z$, we have
\[ \Omega^{nj+i}(M)|\Omega^j(M)^{\otimes n},\]
and since at least one of $\Omega^{nj+i}(M)$ and $\Omega^j(M)$ is non-algebraic, we see that some tensor power of $\Omega^j(M)$ contains a non-algebraic summand; hence $\Omega^j(M)$ is non-algebraic, as required.

Similarly, if $\Omega^i(M^*)\cong \Omega^{-i}(M)^*$ is a summand of $M^{\otimes n}$, then
\[ \Omega^{nj+i}(M^*)|\Omega^j(M)^{\otimes n},\]
and since $\Omega^{nj+i}(M^*)\cong \Omega^{-(nj+i)}(M)^*$, at least one of $\Omega^j(M)$ and $\Omega^{nj+i}(M^*)$ is non-algebraic, and so $\Omega^j(M)$ is non-algebraic.
\end{pf}

\section{The Auslander--Reiten Quiver}

For the basic properties of complexity, we refer to \cite[Proposition 2.2.24]{benson}. One important property that we will use is that the complexities of every module on a particular component of the (stable) Auslander--Reiten quiver is the same. If $B$ is a wild block, then by a theorem of Karin Erdmann in \cite{erdmann1995}, any component $\Gamma$ of $\Gamma_s(B)$ has tree class $A_\infty$. This will be essential in what is to follow.

To prove Theorem \ref{vertexcpcp}, we first introduce the concept of an interlaced component of $\Gamma_s(KG)$. If $\Gamma$ is a component and $\Gamma$ consists either of non-periodic modules or of modules of even periodicity, then for each $M$ in $\Gamma$, the module $\Omega(M)$ does not lie on $\Gamma$. An \emph{interlaced component} is the union of the component $\Gamma$ and the component consisting of the Heller translates of the modules on $\Gamma$. The reason for the name will become clear in the next paragraph.

We begin by co-ordinatizing a non-periodic, interlaced component of $\Gamma_s(KG)$ of type $A_\infty$, which will help immensely in this section. We co-ordinatize according to the following diagram.

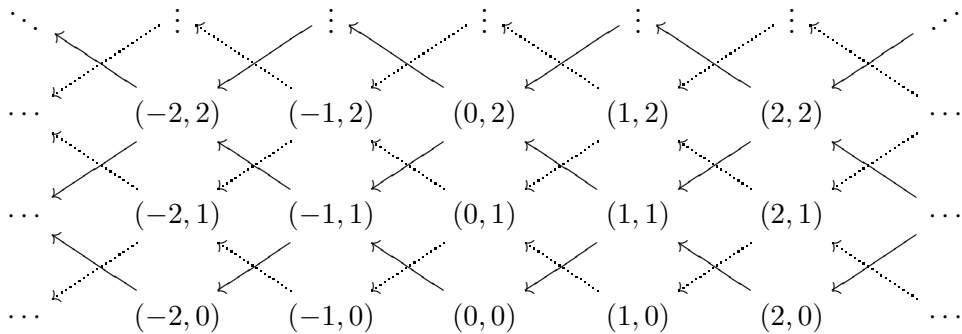
\[\begin{diagram}\node{\ddots}\node{\vdots}\arrow{sw,..}\node{\vdots}\arrow{sw}\node{\vdots}\arrow{sw,..}\node{\vdots}\arrow{sw}\node{\vdots}\arrow{sw,..}\node{\bdots}\arrow{sw}
\\ \node{\cdots}\node{(-2,2)}\arrow{nw}\arrow{sw}\node{(-1,2)}\arrow{nw,..}\arrow{sw,..}\node{(0,2)}\arrow{nw}\arrow{sw}\node{(1,2)}\arrow{nw,..}\arrow{sw,..}\node{(2,2)}\arrow{nw}\arrow{sw}\node{\cdots}\arrow{nw,..}\arrow{sw,..}
\\ \node{\cdots}\node{(-2,1)}\arrow{nw,..}\arrow{sw,..}\node{(-1,1)}\arrow{nw}\arrow{sw}\node{(0,1)}\arrow{nw,..}\arrow{sw,..}\node{(1,1)}\arrow{nw}\arrow{sw}\node{(2,1)}\arrow{nw,..}\arrow{sw,..}\node{\cdots}\arrow{nw}\arrow{sw}
\\ \node{\cdots}\node{(-2,0)}\arrow{nw}\node{(-1,0)}\arrow{nw,..}\node{(0,0)}\arrow{nw}\node{(1,0)}\arrow{nw,..}\node{(2,0)}\arrow{nw}\node{\cdots}\arrow{nw,..}\end{diagram}\]
[Note that this quiver consists of interlaced `diamonds'; when we refer to a diamond of an interlaced component, we mean such a collection of four vertices.]

For the rest of this section, $\Gamma$ will denote an interlaced component of $\Gamma_s(KG)$. Write $M_{(i,j)}$ for the indecomposable module in the $(i,j)$ position on $\Gamma$. (Of course, while $j$ is determined, there is choice over which position on $\Gamma$ is $(0,0)$; we will assume that such a choice is made.)

We recall the following easy result.

\begin{lem}[{{\cite[Proposition 4.12.10]{bensonvol1}}}]\label{splitrest} Let $M$ be an indecomposable module with vertex $Q$, and suppose that $H$ is a subgroup of $G$ not containing any conjugate of $Q$. Then the Auslander--Reiten sequence terminating in $M$ splits upon restriction to $H$.\end{lem}

Notice that, for our interlaced component $\Gamma$ and modules $M_{(i,j)}$, this result becomes the statement that if $H$ does not contain a vertex of $M_{(i,j)}$, then for $i>0$,
\[ M_{(i-1,j)}\res H\oplus M_{(i+1,j)}\res H\cong M_{(i,j+1)}\res H\oplus M_{(i,j-1)}\res H.\]
In particular, this implies that if the modules attached to three of the four vertices in a diamond of $\Gamma$ have known restrictions to $H$, the fourth is uniquely determined.

We also need a slight extension to the result that the complexity of every module on the same component is the same.

\begin{lem}\label{compres} Let $\Gamma$ be an interlaced component of the Auslander--Reiten quiver, and let $H$ be a subgroup of $G$. Then for all $M$ on $\Gamma$, the complexity of $M\res H$ is the same.
\end{lem}
\begin{pf} Let $M$ be a module on $\Gamma$ such that $M\res H$ has the smallest complexity, say $n$. Let
\[ 0\to \Omega^2(M)\to N\to M\to 0\]
be the almost-split sequence terminating in $M$. Restricting this sequence to $H$ yields a short exact sequence whose terms are $KH$-modules. Since $\cx(M\res H)=\cx(\Omega(M)\res H)$, and for any short exact sequence the largest two complexities of the terms are equal, the complexity of $N\res H$ is equal to that of $M\res H$, by minimal choice of $M$. Thus if $L$ is connected to any $\Omega^i(M)$, then $\cx(L\res H)=n$. This holds for any module $M$ such that $\cx(M\res H)=n$, so the restrictions of all modules on the component of $\Gamma_s(KG)$ containing $M$ have the same complexity.\end{pf}

This can be used to prove the next theorem, which is Theorem \ref{vertexcpcp} from the introduction.

\begin{thm}\label{propforthm2} Let $G$ be a finite group and let $\Gamma$ be an interlaced component of $\Gamma_s(KG)$. Suppose that $P$ is a $p$-subgroup such that $P$ does not contain a vertex of any module on $\Gamma$, and that for some $M$ on $\Gamma$, the restriction of $M$ to $P$ is non-periodic. Then $\Gamma$ contains at most one algebraic module and such a module lies at the end of $\Gamma$; i.e., it is $M_{(i,0)}$ for some $i\in \Z$.
\end{thm}
\begin{pf} Since $P$ does not contain a vertex of any module on $\Gamma$, any almost-split sequence involving terms on $\Gamma$ splits upon restriction to $P$, and so we consider all modules $M_{(i,j)}$ to be restricted to $P$. For a co-ordinate $(i,j)$ on $\Gamma$, we attach a collection $[a_1,\dots,a_n]$, which are the non-periodic summands of $M_{(i,j)}\res P$ in a decomposition of $M_{(i,j)}\res P$ into indecomposable summands. We call this collection the \emph{signature} of the vertex $(i,j)$. Notice that, since $M_{(i,j)}=\Omega(M_{(i-1,j)})$, to know the signatures of all vertices in a row, it suffices to know the signature of one of them.

Note that by Lemma \ref{compres}, all modules on $\Gamma$ have non-periodic restriction to $P$, and so the signature of each co-ordinate is non-empty.

By the remarks after Lemma \ref{splitrest}, if we know all signatures of the bottom two rows, we can uniquely determine all signatures of higher rows, since three of the four vertices on each diamond will have known signatures. Also, the second row can be determined from the first row, because of the fact that the signature of $(i,1)$ is equal to the sum of the signatures of $(i-1,0)$ and $(i+1,0)$.

In order to easily express the signatures of the vertices, if $x$ is an element of a signature, then denote by $x^i$ the $i$th Heller translate of $x$. Let $[x_1,\dots,x_t]$ denote the signature of the vertex $(0,0)$. Then the signature of $(i,0)$ is $[x_1^i,x_2^i,\dots,x_t^i]$, and the signature of $(i,1)$ is
\[ [x_1^{i-1},x_2^{i-1},\dots,x_t^{i-1},x_1^{i+1},x_2^{i+1},\dots,x_t^{i+1}],\]
since the almost-split sequence terminating in $M_{(i,0)}$ splits on restriction to $P$. Write $X^i$ for the signature 
\[ [x_1^i,x_2^i,\dots,x_t^i],\]
and write $X^A$ for the signature
\[ \bigcup_{a\in A} X^a.\]

We claim that the signature of $(i,j)$ is
\[ X^{\{i+j,i+j-2,\dots,i-j+2,i-j\}}.\]
To prove this, we firstly note that for $j=0$ and $j=1$ this formula holds. Since we know that the signatures of all vertices are uniquely determined by the first two rows, we simply have to show that it obeys the rule that, for each diamond, the sum of the signatures of the top and bottom vertices equal the sum of the signatures of the left and right vertices. This is true, as the top and bottom vertices' signatures are
\[ X^{\{i+(j+1),i+(j+1)-2,\dots,i-(j+1)+2,i-(j+1)\}}\cup X^{\{i+(j-1),i+(j-1)-2,\dots,i-(j-1)+2,i-(j-1)\}},\]
and the left and right vertices' signatures are
\[ X^{\{(i+1)+j,(i+1)+j-2,\dots,(i+1)-j+2,(i+1)-j\}}\cup X^{\{(i-1)+j,(i-1)+j-2,\dots,(i-1)-j+2,(i-1)-j\}}.\]
These are easily seen to be the same, and so the above formula gives the signature of the vertex $(i,j)$.

If $j\neq 0$, then the signature contains $X^{\{i+j,i+j-2\}}$ and since it cannot be that both $x_1^{i+j}$ and $x_1^{i+j-2}$ are algebraic, $M_{i,j}$ is non-algebraic for all $j>0$.

Finally, at most one of the modules $M_{(i,0)}$ can be algebraic, and so the theorem is proved.\end{pf}

Theorem \ref{propforthm2} can be used to produce results such as Theorem \ref{comp3} and Corollary \ref{p'dim}. In the first case, if the complexity of a module $M$ is at least 3, then the vertex $P$ of $M$ is of $p$-rank 3. By the Alperin--Evens theorem \cite{alperinevens1981}, there is an elementary abelian subgroup $Q$ of $P$ such that $M\res Q$ has complexity 3, and hence for any subgroup $R$ of $Q$ of index $p$, the module $M\res R$ is non-periodic, yielding an appropriate subgroup.

To prove Corollary \ref{p'dim}, recall that a $p'$-dimensional module has a Sylow $p$-subgroup $P$ as a vertex. If $P$ has $p$-rank at least 3, then the result is true by Theorem \ref{comp3}, so $G$ has $p$-rank 2. Let $M$ denote a module on $\Gamma$. By the Alperin--Evens theorem, there is a subgroup $Q$ of $P$ isomorphic with $C_p\times C_p$, such that the complexity of $M\res Q$ is 2, and so $Q$ is a subgroup that satisfies the conditions of Theorem \ref{propforthm2}.

In general it appears difficult to prove a corresponding theorem to Theorem \ref{comp3} for arbitrary $A_\infty$-components of complexity 2. Theorem \ref{propforthm2} places significant restrictions on a possible counterexample to the statement that no non-periodic component of $\Gamma_s(KG)$ from a wild block contains more than one algebraic module.

For dihedral $2$-groups, a similar result -- that a non-periodic component of $\Gamma_s(KG)$ contains at most one algebraic module -- is also true \cite{craventhesis}. The author has not yet considered the semidihedral groups.

\section{Relating Algebraicity and Periodicity}

The results above tell us nothing about the indecomposable modules for $C_p\times C_p$. In this case, there is a very strong conjecture regarding the relationship between algebraic modules and periodic modules, as given in the introduction. We will discuss the computational evidence gathered by the author. We firstly note that neither direction of Conjecture \ref{cpcpper} is obvious, or indeed even known.

The author has constructed all indecomposable modules for $C_3\times C_3$ of dimensions 3 and 6 over $\GF(3)$, and for each of them, has analyzed whether it is algebraic. There are twelve such indecomposable modules of dimension $3$, and over two-hundred absolutely indecomposable modules of dimension $6$. The periodic modules are proved to be algebraic simply by decomposing tensor powers of them. (This incidentally provides hundreds more examples of periodic, algebraic modules.) The non-periodic indecomposable modules can each be proved to be non-algebraic by Corollary \ref{easynonalg}. This fact might be of interest, since it might offer a method by which one half of Conjecture \ref{cpcpper} could be proved.

The author has also constructed all of the 5-dimensional modules for $C_5\times C_5$ over $\GF(5)$, and is in the process of verifying this conjecture for these modules.

All told, thousands of periodic modules for $C_p\times C_p$ are known to be algebraic, and as well as the infinitude of non-algebraic, non-periodic modules provided for by Theorem \ref{algheller}, thousands more low-dimensional non-algebraic, non-periodic modules have been found. (Of course, the non-periodic modules arising from Theorem \ref{algheller} have large dimension in general.)

Moving away from the group $C_p\times C_p$ to general groups, if $G$ is a generalized quaternion group, then $G$ possesses periodic, non-algebraic modules. This is therefore true for any $2$-group with a quaternion subgroup, such as the semidihedral groups. For odd $p$, however, no periodic, non-algebraic indecomposable modules are known. It would be interesting to find such a module, if one exists.

\bibliography{references}

\providecommand{\bysame}{\leavevmode\hbox to3em{\hrulefill}\thinspace}
\providecommand{\MR}{\relax\ifhmode\unskip\space\fi MR }
\providecommand{\MRhref}[2]{%
  \href{http://www.ams.org/mathscinet-getitem?mr=#1}{#2}
}
\providecommand{\href}[2]{#2}
\begin{thebibliography}{10}

\bibitem{alperin1976b}
Jonathan Alperin, \emph{On modules for the linear fractional groups},
  International Symposium on the Theory of Finite Groups, 1974, Tokyo (1976),
  157--163.

\bibitem{alperinevens1981}
Jonathan Alperin and Leonard Evens, \emph{Representations, resolutions and
  {Q}uillen's dimension theorem}, J.~Pure Appl.~Algebra \textbf{22} (1981),
  1--9.

\bibitem{auslandercarlson1986}
Maurice Auslander and Jon Carlson, \emph{Almost-split sequences and group
  rings}, J.~Algebra \textbf{103} (1986), 122--140.

\bibitem{benson}
David Benson, \emph{Modular representation theory: New trends and methods},
  Springer--Verlag, 1984.

\bibitem{bensonvol1}
\bysame, \emph{Representations and cohomology, {I}. {B}asic representation
  theory of finite groups and associative algebras}, Cambridge, 1998.

\bibitem{bensoncarlson1986}
David Benson and Jon Carlson, \emph{Nilpotent elements in the {G}reen ring},
  J.~Algebra \textbf{104} (1986), 329--350.

\bibitem{craventhesis}
David~A.\ Craven, \emph{Algebraic modules for finite groups}, Ph.D. thesis,
  University of Oxford, 2008.

\bibitem{erdmann1995}
Karin Erdmann, \emph{On {A}uslander--{R}eiten components for group algebras},
  J.~Pure Appl.~Algebra \textbf{104} (1995), 149--160.

\bibitem{feit1980}
Walter Feit, \emph{Irreducible modules of $p$-solvable groups},
  Proc.~Sympos.~Pure Math (Santa Cruz, 1979) \textbf{37} (1980), 405--412.

\bibitem{feit}
\bysame, \emph{The representation theory of finite groups}, North--Holland,
  1982.

\end{thebibliography}

\end{document}